# Some Remarks on the $l_1$-Robust Solution of LexRank Problem

**Anna Timonina-Farkas**

**Abstract** Graph-based ranking methods, such as LexRank, are fundamental in Natural Language Processing (NLP) applications like text summarization, as they measure the relative importance of textual units. Building on recent advances in ranking methods for growing and dynamic graphs, we develop a robust variant of LexRank that operates on stochastic similarity graphs with uncertain and expanding structure. Our approach introduces a novel $l_1$-based formulation that captures ambiguity in both transition probabilities and graph size, while maintaining sparsity. The resulting non-convex problem is upper-bounded by a linear program, providing a tractable and interpretable approximation.

## 1 Introduction

In this article, we introduce a robust graph-based method for computing relative importance of textual units for Natural Language Processing. We consider two sets of textual units, i.e., sentences: the set $\mathscr{A}$ contains $N_{\mathscr{A}}$ verified sentences from documents (e.g., from academic literature), while set $\mathscr{B}$ contains $N_{\mathscr{B}}$ generated sentences, which require comparative assessment (i.e., ranking) with respect to the set $\mathscr{A}$. Let $N = N_{\mathscr{A}} + N_{\mathscr{B}}$ and note that the similarity between two sentences can be measured by the cosine between vectors representing these sentences as shown in [2]. Using lower thresholds for similarity and constructing adjacency matrices given these thresholds, one can easily obtain a column-stochastic non-negative transition matrix $P \in \mathbb{R}^{N \times N}$ from the adjacency matrix of the similarity graph by dividing each element by the corresponding row sum (see [2] and Figure 1 for schematic representation).

Anna Timonina-Farkas, PhD
École Polytechnique Fédérale de Lausanne (TOM, EPFL)
EPFL-CDM-MTEI-TOM, Station 5, CH-1015 Lausanne, Switzerland
Institute for Management Development (IMD)
Chem. de Bellerive 23, CH-1007 Lausanne, Switzerland
E-mail: anna.farkas@imd.org, anna.farkas@epfl.ch

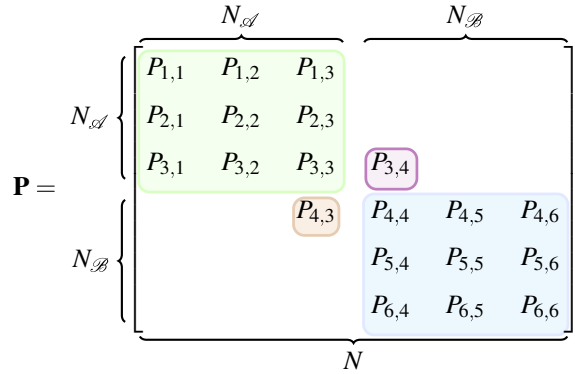

Fig. 1: Transition matrix $P$ encompassing sentences from the academic literature ($N_{\mathscr{A}}$) and the generated sentences ($N_{\mathscr{B}}$). Values $P_{3,4}$ and $P_{4,3}$ represent similarities between generated sentences and verified sentences from academic literature.

Since every sentence is similar at least to itself, all row sums are nonzero. The transportation matrix $P \in \mathbb{R}^{N \times N}$ satisfies $P_{ij} \geq 0$, $\forall i, j = 1, ..., N$ and $\sum_{i=1}^{N} P_{ij} = 1$, $\forall j = 1, ..., N$ by definition. The well-known Perron-Frobenius theorem [3] states that there exists *dominant* (or *principal*) eigenvector $\bar{x} \in \Sigma_N$, where we denote by $\Sigma_N = \{v \in \mathbb{R}^N : \sum_{i=1}^{N} v_i = 1, \ v_i \geq 0\}$ the standard simplex on $\mathbb{R}^N$, i.e., $P\bar{x} = \bar{x}$. Each element of the eigenvector $\bar{x}$ represents the importance (or, ranks) of the corresponding sentence. Looking ahead, we note that the ranks of "true" sentences can be assumed to be high, while the ranks of generated sentences are to be found with respect to the set $\mathscr{A}$.

It is important to note that the number of generated sentences is expected to grow rapidly. One needs to be able to rank these sentences on the fly. However, the dominant vector $\bar{x}$ is not robust and may be highly vulnerable to small changes in matrix $P$. Focusing on Internet-ranking applications, A. Juditsky and B. Polyak in their work [4] formulated a robust optimization problem, allowing matrix $P$ to vary according to the law $P + \xi$ under some conditions on $\xi$. For this matrix, they found the stationary distribution $\bar{x}$ robust to variations in Internet links and, therefore, stable with respect to small changes of $P$. However, in their work, the size of matrix $\xi$ was assumed to be the same as of matrix $P$, i.e., $N \times N$, meaning that growing in the number of nodes was not considered in the network (further, this formulation is referred to as *a fixed-size model*). Therefore, A. Timonina-Farkas and R. Seifert [5] made a step forward allowing the network to have a fixed growth rate via considering the following uncertain transition matrix:

$$Q = \begin{pmatrix} P + \xi & \zeta \\ \psi & \chi \end{pmatrix}, \tag{1}$$

where $P$ is the column-stochastic transportation matrix describing the current state of the network with $N$ sentences; $\xi$ is the matrix describing variations in links of the initial network; $\zeta$, $\psi$ and $\chi$ are matrices describing similarities with $M$ new sentences, which may appear in the future ($\xi$ is of the size $N \times N$, $\psi$ is of the size $M \times N$, $\zeta$ is of the size $N \times M$ and $\chi$ is of the size $M \times M$).

Note that in LexRank application, changes of matrix $P$ can happen solely due to the changes in the number of nodes because the similarities between initial $N$ sentences do not change. Thus, in this work we follow the approach of A. Timonina-Farkas and R. Seifert [5] and propose an $l_1$-norm formulation for robust ranking. We show that the resulting problem is independent of the growth rate given a particular ambiguity set.

As matrices $P$ and $Q$ must be column-stochastic, $\xi$, $\zeta$, $\psi$ and $\chi$ must satisfy the following properties:

$$\begin{cases} \xi_{ij} \geq -P_{ij}, \ \forall i,j = 1,...,N; \\ \psi_{ij} \geq 0, \ \forall i = 1,...,M, \ j = 1,...,N; \\ \zeta_{ij} \geq 0, \ \forall i = 1,...,N, \ j = 1,...,M; \\ \chi_{ij} \geq 0, \ \forall i,j = 1,...,M. \end{cases} \quad (2)$$

Also, the following properties must hold:

$$\begin{cases} \sum_{i=1}^{N} \xi_{ij} + \sum_{i=1}^{M} \psi_{ij} = 0, \ \forall j = 1,...,N \\ \sum_{i=1}^{N} \zeta_{ij} + \sum_{i=1}^{M} \chi_{ij} = 1, \ \forall j = 1,...,M, \end{cases} \quad (3)$$

saying that every column of the matrix $Q$ sums up to 1 (notice, that here we use the fact that $P$ is also column-stochastic and $\sum_{i=1}^{N} P_{ij} = 1, \ \forall j = 1,...,N$).

The function $\max_{Q \in \Xi} ||Qx - x||_{(*)}$, where $|| \cdot ||_{(*)}$ is some norm, can be seen as a measure of "goodness" of a vector $x$ as a common dominant eigenvector of the family $\Xi$, where $\Xi$ stands for the set of perturbed stochastic matrices of the form (1) under conditions (2) and (3). In the work of A. Timonina-Farkas and R. Seifert [5], the authors considered a fixed growth rate and $l_2$- and Frobenius-norm formulations. In this work, we consider an $l_1$-norm formulation, which does not discourage sparsity in transition probabilities. Indeed, consider the uncertainty matrix $\psi$, which describes links from $N$ existing sentences to $M$ new ones. All elements of the matrix $\psi$ are non-negative. If we reduce one positive element from $j$-th column of this matrix by a small enough $\delta$, the norm $||[\psi]_j||_1$ and the sum $\sum_{i,j} |\psi_{ij}|$ decrease by this $\delta$, regardless of the value of the element we decrease. This means, that $l_1$-norm formulation (5) does not make a preference which transition probabilities to decrease in order to satisfy the constraints. By this, a lot of transition probabilities of the matrix $Q$ can result in being zeros. In contrast, for $l_2$- and Frobenius-norm formulations, the reduction of larger terms of the matrix $\psi$ by $\delta$ results in a much greater reduction in norms than doing so with smaller terms. Therefore, $l_2$- and Frobenius-norm formulations discourage sparsity by yielding diminishing reductions for elements closer to zero.

Let us denote by $x = \begin{pmatrix} x^{(1)} \\ x^{(2)} \end{pmatrix}$ a feasible point, which is a candidate for the common dominant eigenvector of the family $\Xi$. Let $x^{(1)}$ be of the size $N \times 1$ and $x^{(2)}$ be of the size $M \times 1$. Hence, $x$ is of the size $(N + M) \times 1$. Notice, that the vector $x$ must belong to the standard simplex $\Sigma_{N+M} = \{v \in \mathbb{R}^{N+M} : \sum_{i=1}^{N+M} v_i = 1, \ v_i \geq 0\}$, that means $x_i^{(1)} \geq 0, \ \forall i = 1,...,N, \ x_j^{(2)} \geq 0, \ \forall j = 1,...,M$ and $\sum_{i=1}^{N} x_i^{(1)} + \sum_{j=1}^{M} x_j^{(2)} = 1$ (i.e., $||x^{(1)}||_1 + ||x^{(2)}||_1 = 1$). We say that the vector $\hat{x}$ is an $l_1$-robust solution of the eigenvector problem on $\Xi$ for any growth rate $M$ if

$$\hat{x} \in \underset{x \in \Sigma_{N+M}}{\mathrm{Argmin}} \left\{ \max_{Q \in \Xi,\ M \in \mathbb{N}} ||Qx - x||_1 \right\}, \tag{4}$$

where $\mathbb{N} = \{0, 1, 2, ...\}$ is the set of non-negative integers.

The reasonable choice of the uncertainty set $\Xi$ would impose some bounds on the column-wise norms of matrices $\xi$, $\zeta$, $\psi$ and $\chi$, meaning that the perturbation in links of the current and future states of the network would be bounded: i.e. $||[\xi]_j|| \leq \varepsilon_j^{(\xi)}$, $||[\psi]_j|| \leq \varepsilon_j^{(\psi)}$, $||[\zeta]_j|| \leq \varepsilon_j^{(\zeta)}$ and $||[\chi]_j|| \leq \varepsilon_j^{(\chi)}$, where $[\cdot]_j$ denotes the $j$-th column of a matrix. Moreover, the *total uncertainty budget* for matrices $\xi$, $\zeta$, $\psi$ and $\chi$ could be fixed (see [4]): this would imply constraints on the overall possible perturbations of the transportation matrix $P$. By solving the optimization problem (4), one would protect the rank vector $\hat{x}$ against high fluctuation in case of link or node perturbations influencing the transportation matrix $P$.

## 2 Problem Formulation and its Bounds

We say that the vector $\hat{x}^{(l_1)}$ is the $l_1$-robust solution of the eigenvector problem on the set of perturbed matrices $\Xi^{(l_1)}$ for any growth rate $M$, if

$$\begin{aligned}
\hat{x}^{(l_1)} \in \underset{x \in \Sigma_{N+M}}{\mathrm{Argmin}} \quad &\left\{ \max_{Q, M} ||Qx - x||_1 : \right. \\
&Q \text{ is column-stochastic, } M \in \{0, 1, 2, ...\} \\
&||[\xi]_j||_1 \leq \varepsilon_j^{(\xi)}, \ \forall j = 1, ..., N, \text{ and } \sum_{i,j} |\xi_{ij}| \leq \varepsilon^{(\xi)}, \\
&||[\psi]_j||_1 \leq \varepsilon_j^{(\psi)}, \ \forall j = 1, ..., N, \text{ and } \sum_{i,j} |\psi_{ij}| \leq \varepsilon^{(\psi)}, \\
&||[\zeta]_j||_1 \leq \varepsilon_j^{(\zeta)}, \ \forall j = 1, ..., M, \text{ and } \sum_{i,j} |\zeta_{ij}| \leq \varepsilon^{(\zeta)}, \\
&\left. ||[\chi]_j||_1 \leq \varepsilon_j^{(\chi)}, \ \forall j = 1, ..., M, \text{ and } \sum_{i,j} |\chi_{ij}| \leq \varepsilon^{(\chi)}, \right\},
\end{aligned} \tag{5}$$

where we have $l_1$-norm constraints on each column $[\cdot]_j$, $\forall j$ of matrices $\xi$, $\zeta$, $\psi$ and $\chi$, as well as we bound the sum of absolute values of all elements of these matrices. We note that $\varepsilon_j^{(\xi)} \leq 1$ and $\varepsilon_j^{(\psi)} \leq 1$ due to the column-stochasticity of the transition matrix. We also note that $\varepsilon_j^{(\zeta)} + \varepsilon_j^{(\chi)} = 1$, $\forall j$ and, thus, $\varepsilon^{(\zeta)} + \varepsilon^{(\chi)} = M$. Further, we assume that the total uncertainty budgets $\varepsilon^{(\xi)}$ and $\varepsilon^{(\psi)}$ (and also $\varepsilon_j^{(\xi)}$ and $\varepsilon_j^{(\psi)}$) are independent of $M$. Oppositely, budgets $\varepsilon^{(\chi)}$ and $\varepsilon^{(\zeta)}$ may increase as $M$ grows. Note that the problem (5) is not convex-concave and can be written as follows:

$$\hat{x}^{(l_1)} \in \underset{x \in \Sigma_{N+M}}{\mathrm{Argmin}} \max_{\substack{Q \in \Xi^{(l_1)} \\ M \in \mathbb{N}}} ||Qx - x||_1.$$

Now, consider the inner optimization problem:

$$\phi_1(x) = \max_{\substack{Q \in \Xi^{(l_1)} \\ M \in \mathbb{N}}} \|Qx - x\|_1.$$

**Proposition 1** *Optimal value $\phi_1(\hat{x}^{(l_1)})$ can be bounded from above by the optimal value of the following convex optimization problem with $\varepsilon_1 = \varepsilon^{(\xi)} + \varepsilon^{(\psi)}$:*

$$\phi_1(\hat{x}^{(l_1)}) \leq \min_{x \in \Sigma_{N+M}} \left\{ \|Px^{(1)} - x^{(1)}\|_1 + \varepsilon_1 \|x^{(1)}\|_{(a)} + \max_{\substack{z \in \mathbb{R}^M, M \in \mathbb{N} \\ \|z\|_1 \leq \varepsilon^{(\zeta)} + \varepsilon^{(\chi)} + M \\ |z_j| \leq \varepsilon_j^{(\zeta)} + \varepsilon_j^{(\chi)} + 1}} z^T x^{(2)} \right\}, \quad (6)$$

*where $\|x^{(1)}\|_{(a)} = \min_{\lambda + \mu = x^{(1)}} \left\{ \|\lambda\|_\infty + \sum_{j=1}^N \frac{\varepsilon_j^{(\xi)} + \varepsilon_j^{(\psi)}}{\varepsilon^{(\xi)} + \varepsilon^{(\psi)}} |\mu_j| \right\}.$*

*Proof* See the Appendix 6.1 for the proof.

Note that $\|x^{(2)}\|_{(b)} = 0$ if there are no new generated sentences in the network (i.e., if $M = 0$). In this case, the optimization problem (6) coincides with the $l_1$-reformulation proposed by A. Juditsky and B. Polyak in the work [4]. Furthermore, in line with the work of A. Timonina-Farkas and R.W. Seifert [5], similar upper bounds hold for $l_2$ and Frobenius norms.

Optimization problem (5) accounts for uncertainties in similarities between current and future (not yet generated) sentences: these uncertainties are incorporated via matrices $\xi$, $\psi$, $\zeta$ and $\chi$, which worst-case realization gives us the opportunity to compute the robust LexRank. This optimization problem differs from robust formulations corresponding to the fixed-size network model proposed by A. Juditsky and B. Polyak in the work [4], who studied uncertainties implied by the matrix $P + \xi$, describing variations in links of existing pages with constant network size.

**Proposition 2** *Optimization problem (5) under conditions (2) and (3) imposes upper bounds on the fixed-size network model in the following sense:*

$$\phi_1(x) = \max_{\substack{Q \in \Xi^{(l_1)} \\ M \in \mathbb{N}}} \|Qx - x\|_1 \geq \max_{\substack{\mathbb{1}_N^T[\xi]_j = 0 \\ \|[\xi]_j\|_1 \leq \varepsilon_j^{(\xi)} \\ \Sigma_{i,j} |\xi_{ij}| \leq \varepsilon^{(\xi)}}} \|(P + \xi)x^{(1)} - x^{(1)}\|_1, \quad (7)$$

*Proof* See the Appendix 6.2 for the proof.

In the next section, we consider the optimization problem (6) and develop methods for its solution.

## 3 Influence of newly generated sentences

The dimensionality of the optimization problem (6) can be reduced, as the problem can be subdivided into two separate smaller-size optimization problems as stated in the next Lemma 1.

**Lemma 1** *The optimal solution* $\hat{x}^{(l_1)} = \begin{pmatrix} \hat{x}^{(1)} \\ \hat{x}^{(2)} \end{pmatrix}$ *of the optimization problem (6) is independent of the growth rate M and is equal to the following:*

$$\hat{x}^{(1)} \in \underset{x^{(1)} \in \Sigma_N}{\text{Argmin}} \left\{ \|Px^{(1)} - x^{(1)}\|_1 + \varepsilon_1 \|x^{(1)}\|_{(a)} \right\}, \ \hat{x}^{(2)} = 0, \tag{8}$$

*where* $\Sigma_N = \{ v \in \mathbb{R}^N : \sum_{i=1}^N v_i, \ v_i \geq 0 \}.$

*Proof* Consider the optimization problem (6). Based on the fact, that $\|x^{(1)}\|_1 + \|x^{(2)}\|_1 = 1$ and $x^{(1)} \geq 0$, $x^{(2)} \geq 0$, let us make the following change of variables with $s \in [0, 1]$:

$$x^{(1)} = sy^{(1)}, \ x^{(2)} = (1-s)y^{(2)}, \tag{9}$$

where $\|y^{(1)}\|_1 = 1$, $\|y^{(2)}\|_1 = 1$ and $y^{(1)} \geq 0$, $y^{(2)} \geq 0$. Further, $\|x^{(1)}\|_1 = s$ and $\|x^{(2)}\|_1 = 1-s$, $\forall s \in [0, 1]$. Note that simplex constraints on $x = \begin{pmatrix} x^{(1)} \\ x^{(2)} \end{pmatrix}$ are satisfied and that the decision about the value of $s$ is due to the outer optimization problem (see (6)). Therefore, the optimization problem (6) can be equivalently rewritten as follows:

$$\min_{s \in [0,1]} \left\{ s \min_{y^{(1)} \in \Sigma_N} \left( \|Py^{(1)} - y^{(1)}\|_1 + \varepsilon_1 \|y^{(1)}\|_{(a)} \right) + (1-s) \min_{y^{(2)} \in \Sigma_M} \max_{\substack{z \in \mathbb{R}^M, M \in \mathbb{N} \\ \|z\|_1 \leq \varepsilon^{(\zeta)} + \varepsilon^{(\chi)} + M \\ |z_j| \leq \varepsilon_j^{(\zeta)} + \varepsilon_j^{(\chi)} + 1}} z^T y^{(2)} \right\}.$$

Next, we show that the following holds $\forall \bar{M}$:

$$\min_{y^{(1)} \in \Sigma_N} \left( \|Py^{(1)} - y^{(1)}\|_1 + \varepsilon_1 \|y^{(1)}\|_{(a)} \right) < \min_{y^{(2)} \in \Sigma_{\bar{M}}} \max_{\substack{z \in \mathbb{R}^{\bar{M}} \\ \|z\|_1 \leq \varepsilon^{(\zeta)} + \varepsilon^{(\chi)} + \bar{M} \\ |z_j| \leq \varepsilon_j^{(\zeta)} + \varepsilon_j^{(\chi)} + 1}} z^T y^{(2)}. \tag{10}$$

If this is the case, then the inequality also holds at optimality. To prove it, we note that the inner optimization in the right-hand side is equivalent to the following:

$$\max_{\substack{z \in \mathbb{R}^{\bar{M}} \\ \|z\|_1 \leq \varepsilon^{(\zeta)} + \varepsilon^{(\chi)} + \bar{M} \\ |z_j| \leq \varepsilon_j^{(\zeta)} + \varepsilon_j^{(\chi)} + 1}} z^T y^{(2)} = (\varepsilon^{(\zeta)} + \varepsilon^{(\chi)} + \bar{M}) \|y^{(2)}\|_{(b)}$$

where $\|y^{(2)}\|_{(b)} = \min_{\lambda + \mu = y^{(2)}} \left\{ \|\lambda\|_\infty + \sum_{j=1}^{\bar{M}} \frac{\varepsilon_j^{(\zeta)} + \varepsilon_j^{(\chi)} + 1}{\varepsilon^{(\zeta)} + \varepsilon^{(\chi)} + \bar{M}} |\mu_j| \right\}$ according to the Lemma 2 in Appendix. Furthermore, as proven in Appendix 6.3 and because $\varepsilon^{(\zeta)} + \varepsilon^{(\chi)} =$

$\bar{M}$ and $\varepsilon_j^{(\zeta)} + \varepsilon_j^{(\chi)} = 1$ due to the column-stochasticity of the transition matrix, the following exact solution holds: $\min_{y^{(2)} \in \Sigma_{\bar{M}}} \|y^{(2)}\|_{(b)} = \frac{1}{\bar{M}}$. Thus, for the inequality (10) to be satisfied, it is sufficient that

$$\min_{y^{(1)} \in \Sigma_N} \left\{ \|Py^{(1)} - y^{(1)}\|_1 + \varepsilon_1 \|y^{(1)}\|_{(a)} \right\} < \frac{\varepsilon^{(\zeta)} + \varepsilon^{(\chi)} + \bar{M}}{\bar{M}} = 2.$$

Note that $\min_{y^{(1)} \in \Sigma_N} \left\{ \|Py^{(1)} - y^{(1)}\|_1 + \varepsilon_1 \|y^{(1)}\|_{(a)} \right\} \leq \sum_{j=1}^N (\varepsilon_j^{(\xi)} + \varepsilon_j^{(\psi)}) |\bar{y}_j^{(1)}|$, where the last inequality holds due to the definition of the norm $\|\cdot\|_{(a)}$ and where $\bar{y}^{(1)}$ is such a vector, that $\bar{y}^{(1)} = P\bar{y}^{(1)}$ and $\bar{y}^{(1)} \in \Sigma_N$ (it exists due to the Perron-Frobenius theorem). Therefore, condition (10) is always satisfied because $\varepsilon_j^{(\xi)} + \varepsilon_j^{(\psi)} \leq 2$ for all $j$ due to matrix stochasticity and $\bar{y}^{(1)} \in \Sigma_N$. Due to this, the optimal $s$ is equal to one and the statement of the Lemma follows. $\qquad\square$

The solution in the Lemma 1 depends neither on the number $M$ of new sentences in the network nor on the uncertainty levels $\varepsilon^{(\zeta)}$ and $\varepsilon^{(\chi)}$, which makes the $l_1-$ formulation convenient for the analysis. Note that it is not necessarily the case of $l_2-$ robust reformulation presented by A. Timonina-Farkas and R.W. Seifert [5]. Indeed, their condition requires that the uncertainty about future sentences $\varepsilon^{(\zeta)} + \varepsilon^{(\chi)}$ grows faster than $(\sqrt{M} - 1)$ as the number $M$ of new pages increases. For the case $M = 0$ and $M = 1$, this condition is automatically satisfied. For higher dimensions (i.e. $M > 1$), this condition is statistically tested from real-world data (see [5]). Similarly, for the Frobenius-norm formulation, the condition requires that the uncertainty about future pages $\varepsilon^{(\zeta)} + \varepsilon^{(\chi)}$ grows faster than $(\varepsilon^{(\xi)} + \varepsilon^{(\psi)})\sqrt{M} - 1$ as the number $M$ of new pages increases.

Further, we focus on the solution of the problem (8), which is formulated in the following general form:

$$\min_{\substack{x \in \Sigma_N \\ \lambda + \mu = x}} \left\{ \|Px - x\|_1 + \varepsilon_1 \|\lambda\|_\infty + \sum_{j=1}^N \varepsilon_j |\mu_j| \right\}. \tag{11}$$

Here, with some abuse of notations, we denote $\varepsilon_j^{(\xi)} + \varepsilon_j^{(\psi)}$ by $\varepsilon_j$ with $\varepsilon = (\varepsilon_1, ..., \varepsilon_N)^T$. Clearly, the optimization (11) can be written in a linear programming form (LP), i.e.,

$$\min_{x, s, t, r, \mu} \quad \mathbb{1}^T q + \varepsilon_1 t + \varepsilon^T r, \tag{12}$$
$$\text{subject to } -s \leq Px - x \leq s, \ s \geq 0$$
$$-t\mathbb{1} \leq x - \mu \leq t\mathbb{1}, \ t \geq 0$$
$$-r \leq \mu \leq r, \ r \geq 0, \ x \geq 0, \ \mathbb{1}^T x = 1.$$

In the next section, we demonstrate the solution of this LP for a LexRank problem.

## 4 Numerical insights

We conduct small- to medium-dimensional tests demonstrating the solution of the optimization problem (11). The high-dimensional analysis is similar to the work of A. Timonina-Farkas and R.W. Seifert [5] and, therefore, is not provided in this article.

### 4.1 Small-dimensional tests

Any cluster of documents can be represented by a cosine similarity matrix where each entry is the similarity score between the corresponding sentence pair. For a small-dimensional test, we use the dataset of 11 sentences provided in the article of E. Günes and D.R. Radev [2]. The dataset is shown in Table 1. Further, intra-sentence cosine similarities between sentences are computed as in E. Günes and D.R. Radev [2] and shown in Figure 2.

Table 1: Sentences from documents with IDs as of E. Günes and D.R. Radev [2].

| $N°$ | ID | Sentence |
|---|---|---|
| 1 | d1s1 | Iraqi Vice President Taha Yassin Ramadan announced today, Sunday, that Iraq refuses to back down from its decision to stop cooperating with disarmament inspectors before its demands are met. |
| 2 | d2s1 | Iraqi Vice President Taha Yassin Ramadan announced today, Thursday, that Iraq rejects cooperating with the United Nations except on the issue of lifting the blockade imposed upon it since the year 1990. |
| 3 | d2s2 | Ramadan told reporters in Baghdad that "Iraq cannot deal positively with whoever represents the Security Council unless there was a clear stance on the issue of lifting the blockade off of it." |
| 4 | d2s3 | Baghdad had decided late last October to completely cease cooperating with the inspectors of the United Nations Special Commission (UNSCOM), in charge of disarming Iraq's weapons, and whose work became very limited since the fifth of August, and announced it will not resume its cooperation with the Commission even if it were subjected to a military operation. |
| 5 | d3s1 | The Russian Foreign Minister, Igor Ivanov, warned today, Wednesday against using force against Iraq, which will destroy, according to him, seven years of difficult diplomatic work and will complicate the regional situation in the area. |
| 6 | d3s2 | Ivanov contended that carrying out air strikes against Iraq, who refuses to cooperate with the United Nations inspectors, "will end the tremendous work achieved by the international group during the past seven years and will complicate the situation in the region." |
| 7 | d3s3 | Nevertheless, Ivanov stressed that Baghdad must resume working with the Special Commission in charge of disarming the Iraqi weapons of mass destruction (UNSCOM). |

| 8 | d4s1 | The Special Representative of the United Nations Secretary-General in Baghdad, Prakash Shah, announced today, Wednesday, after meeting with the Iraqi Deputy Prime Minister Tariq Aziz, that Iraq refuses to back down from its decision to cut off cooperation with the disarmament inspectors. |
|---|------|---|
| 9 | d5s1 | British Prime Minister Tony Blair said today, Sunday, that the crisis between the international community and Iraq "did not end" and that Britain is still "ready, prepared, and able to strike Iraq." |
| 10 | d5s2 | In a gathering with the press held at the Prime Minister's office, Blair contended that the crisis with Iraq "will not end until Iraq has absolutely and unconditionally respected its commitments" towards the United Nations. |
| 11 | d5s3 | A spokesman for Tony Blair had indicated that the British Prime Minister gave permission to British Air Force Tornado planes stationed in Kuwait to join the aerial bombardment against Iraq. |

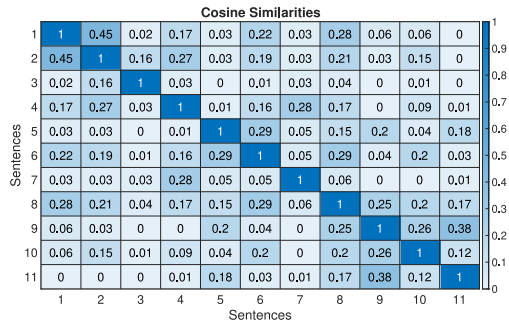

Fig. 2: Heatmap for intra-sentence cosine similarities between sentences in Table 1.

Given the thresholds of 0.1, 0.2, and 0.3 as lower bounds for sentence similarities, we construct the adjacency matrix and, then, compute the transition matrix $P$ shown in Figure 6.

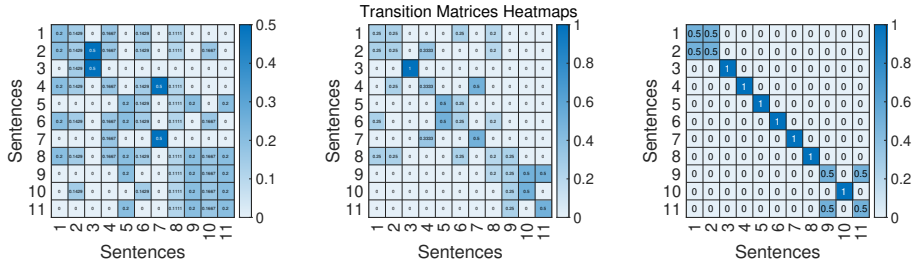

Fig. 3: Heatmap for transition matrices given 0.1, 0.2, and 0.3 lower bounds for sentence similarities.

LexRank scores for the graphs in Figure 6 and their robust counterparts for different values of $\varepsilon_1$ and $\varepsilon_j$ are provided in Tables 2-4. Note that the values are normalized so that the largest value of each column is equal to one.

Table 2: LexRank with and without robustness, $\varepsilon_1 = 0.01$ and $\varepsilon_j = 0.01$, $\forall j$.

| IDs | LR (0.1) | Robust LR (0.1) | LR (0.2) | Robust LR (0.2) | LR (0.3) | Robust LR (0.3) |
|-----|----------|-----------------|----------|-----------------|----------|-----------------|
| d1s1 | 0.6007 | 0.5556 | 0.6944 | 0.8000 | 1 | 1 |
| d2s1 | 0.8466 | 0.7778 | 0.7317 | 0.8000 | 1 | 1 |
| d2s2 | 0.3491 | 0.2222 | 0.6773 | 1.0000 | 1 | 1 |
| d2s3 | 0.7520 | 0.6667 | 0.6550 | 0.6667 | 1 | 1 |
| d3s1 | 0.5907 | 0.5556 | 0.4344 | 0.4000 | 1 | 1 |
| d3s2 | 0.7993 | 0.7778 | 0.8718 | 0.8000 | 1 | 1 |
| d3s3 | 0.3548 | 0.2222 | 0.4993 | 0.4000 | 1 | 1 |
| d4s1 | 1.0000 | 1.0000 | 1.0000 | 1.0000 | 1 | 1 |
| d5s1 | 0.5921 | 0.5556 | 0.7399 | 0.8000 | 1 | 1 |
| d5s2 | 0.6910 | 0.6667 | 0.6967 | 0.4000 | 1 | 1 |
| d5s3 | 0.5921 | 0.5556 | 0.4501 | 0.4000 | 1 | 1 |

Table 3: LexRank with and without robustness, $\varepsilon_1 = 5$ and $\varepsilon_j = 5$, $\forall j$.

| IDs | LR (0.1) | Robust LR (0.1) | LR (0.2) | Robust LR (0.2) | LR (0.3) | Robust LR (0.3) |
|-----|----------|-----------------|----------|-----------------|----------|-----------------|
| d1s1 | 0.60070 | 0.70437 | 0.69440 | 0.93333 | 1 | 1 |
| d2s1 | 0.84660 | 1.00000 | 0.73170 | 1.00000 | 1 | 1 |
| d2s2 | 0.34910 | 0.28571 | 0.67730 | 1.00000 | 1 | 1 |
| d2s3 | 0.75200 | 1.00000 | 0.65500 | 1.00000 | 1 | 1 |
| d3s1 | 0.59070 | 0.81746 | 0.43440 | 0.63333 | 1 | 1 |
| d3s2 | 0.79930 | 1.00000 | 0.87180 | 1.00000 | 1 | 1 |
| d3s3 | 0.35480 | 0.59127 | 0.49930 | 0.83333 | 1 | 1 |
| d4s1 | 1.00000 | 1.00000 | 1.00000 | 1.00000 | 1 | 1 |
| d5s1 | 0.59210 | 1.00000 | 0.73990 | 1.00000 | 1 | 1 |
| d5s2 | 0.69100 | 0.95619 | 0.69670 | 0.55000 | 1 | 1 |
| d5s3 | 0.59210 | 1.00000 | 0.45010 | 0.55000 | 1 | 1 |

Table 4: LexRank with and without robustness, $\varepsilon_1 = 10$ and $\varepsilon_j = 10$, $\forall j$.

| IDs | LR (0.1) | Robust LR (0.1) | LR (0.2) | Robust LR (0.2) | LR (0.3) | Robust LR (0.3) |
|-----|----------|-----------------|----------|-----------------|----------|-----------------|
| d1s1 | 0.60070 | 0.70437 | 0.69440 | 0.93333 | 1 | 1 |
| d2s1 | 0.84660 | 1.00000 | 0.73170 | 1.00000 | 1 | 1 |
| d2s2 | 0.34910 | 1.00000 | 0.67730 | 1.00000 | 1 | 1 |
| d2s3 | 0.75200 | 1.00000 | 0.65500 | 1.00000 | 1 | 1 |
| d3s1 | 0.59070 | 1.00000 | 0.43440 | 1.00000 | 1 | 1 |
| d3s2 | 0.79930 | 1.00000 | 0.87180 | 1.00000 | 1 | 1 |
| d3s3 | 0.35480 | 1.00000 | 0.49930 | 1.00000 | 1 | 1 |
| d4s1 | 1.00000 | 1.00000 | 1.00000 | 1.00000 | 1 | 1 |
| d5s1 | 0.59210 | 1.00000 | 0.73990 | 1.00000 | 1 | 1 |
| d5s2 | 0.69100 | 1.00000 | 0.69670 | 1.00000 | 1 | 1 |
| d5s3 | 0.59210 | 1.00000 | 0.45010 | 1.00000 | 1 | 1 |

### 4.2 Medium-dimensional tests

Finally, beside testing the method on a small dataset from the article of E. Günes and D.R. Radev [2], we also test the method by automatically extracting 1032 sentences from the article of A. Timonina-Farkas and R. W. Seifert [5] and ranking these sentences (see Figure 4). Three transition matrices are constructed using the thresholds of 0.05, 0.15, and 0.2 for the cosine similarity. As expected, a larger ambiguity set leads to a greater fraction of equalized ranks (see Figure 4).

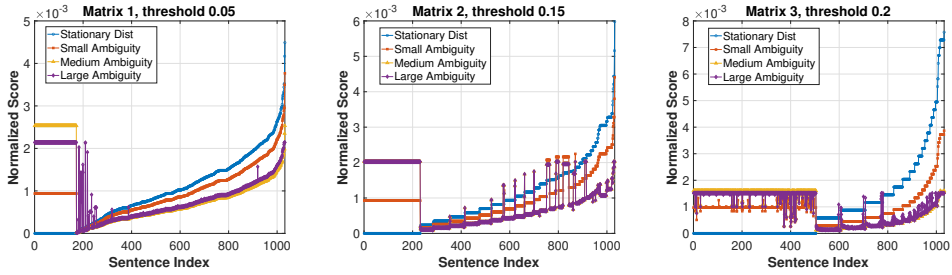

Fig. 4: LexRank with robustness and matrix $P$: 1) no ambiguity $\varepsilon_1 = 0$ and $\varepsilon_j = 0$, $\forall j$; 2) medium ambiguity $\varepsilon_1 = 2$ and $\varepsilon_j = 2$, $\forall j$; 3) large ambiguity $\varepsilon_1 = 7$ and $\varepsilon_j = 7$, $\forall j$.

Next, we show the performance of LexRank algorithms in in the presence of AI-generated sentences. For this, in addition to $N_{\mathscr{A}} = 1032$ sentences from academic literature, we generate $N_{\mathscr{B}} = 2000$ random sentences based on the templates:

Table 5: Sentence templates about LexRank.

| $N°$ | Sentence |
|---|---|
| 1 | LexRank is a graph-based method for text summarization. |
| 2 | The LexRank algorithm uses sentence similarity to build a graph. |
| 3 | In LexRank, nodes represent sentences and edges reflect cosine similarity. |
| 4 | LexRank is related to the PageRank algorithm from web search. |
| 5 | An advantage of LexRank is its ability to find central sentences in a document. |
| 6 | LexRank applies eigenvector centrality to identify important sentences. |
| 7 | One drawback of LexRank is sensitivity to sentence similarity thresholds. |
| 8 | LexRank can be used for extractive summarization in news articles. |
| 9 | Researchers often compare LexRank with TextRank and other summarization algorithms. |
| 10 | LexRank is unsupervised, meaning it does not require labeled training data. |

We observe that the LexRank algorithm, which is based solely on the dominant eigenvector computation (i.e., stationary distribution), may assign high weights to the generated sentences as it is seen in Figure 5. The ranks are strongly dependent on the chosen similarity threshold. For lower thresholds, generated sentences can be weighted higher in comparison to sentences from academic literature. For higher thresholds, the algorithm detects generated sentences by assigning them zero weights. This, however, happens at a cost of the order of ranks in the set $\mathscr{A}$ and implies the necessity to carefully choose the threshold for the similarity.

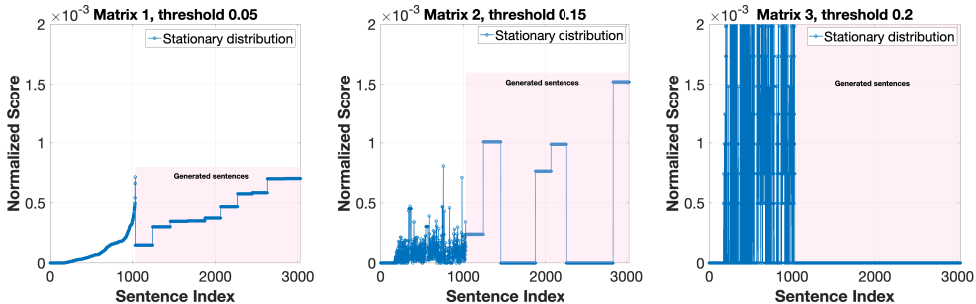

Fig. 5: LexRank for sentences from academic research and generated sentences.

If the primary objective is not the relative ranking of sentences in set $\mathscr{A}$, but rather a comparative assessment of those in set $\mathscr{B}$, our optimization approach can be applied by leveraging the knowledge that the first $N_{\mathscr{A}} = 1032$ sentences originate from academic literature. This leads to the following optimization problem, whose solution $\bar{x} = \frac{x}{\sum_i x}$ provides a feasible point and, consequently, an upper bound for problem (12).

$$\min_{x,s,t,r,\mu} \quad \mathbb{1}^T q + \varepsilon_1 t + \varepsilon^T r, \tag{13}$$

$$\text{subject to} -s \leq Px - x \leq s, \; s \geq 0$$

$$-t\mathbb{1} \leq x - \mu \leq t\mathbb{1}, \; t \geq 0, \; -r \leq \mu \leq r, \; r \geq 0, \; x \geq 0,$$

$$x_1 = 1, \dots x_{N_{\mathscr{A}}} = 1, \; x_{N_{\mathscr{A}}+1} \leq 1, \dots, x_{N_{\mathscr{A}}+N_{\mathscr{B}}} \leq 1.$$

In this setting, the ranks of sentences in $\mathscr{B}$ are evaluated relative to the academic literature and remain less than or equal to those in set $\mathscr{A}$ (see Figure 6). High-quality generated sentences are those whose ranks approach those of the academic sentences. An important direction for future research is to examine whether these comparative ranks can serve as predictive features (alongside other textual metrics) for distinguishing AI-generated from human-written sentences.

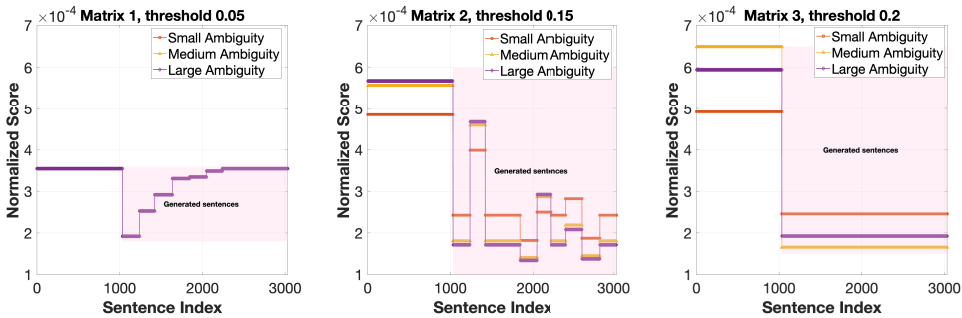

Fig. 6: Comparative LexRank for sentences from academic research and generated sentences.

## 5 Conclusion

In this work, we develop an $l_1$-robust formulation of LexRank, explicitly accounting for uncertainties in both existing and newly generated sentences. We show that, under the $l_1$-norm, the optimal robust solution for existing sentences is independent of the number of newly generated sentences ($M$) and the uncertainty levels associated with them. This property substantially simplifies the analysis and computation of robust rankings. We provide a bound linking our growing-network model to the classical fixed-size robust LexRank formulation, illustrating how the $l_1$-robust approach generalizes existing results while retaining tractability. The robust optimization problem can be efficiently solved via linear programming, and numerical experiments on small- and medium-sized datasets demonstrated the effect of robustness: as the uncertainty budgets increase, rankings become more conservative, reducing sensitivity to small perturbations in the transition matrix. Our results highlight that the $l_1$-robust LexRank formulation provides a practical and scalable approach for extracting stable sentence rankings in networks subject to uncertainty, offering a foundation for further studies in robust text summarization and ranking under network growth.

# References

1. El Ghaoui, L., Lebret, H. *Robust Solutions to Least-Squares Problems with Uncertain Data.* SIAM Journal on Matrix Analysis and Applications, Volume 18(4), pp. 1035–1064 (1997).
2. Günes, E. Radev, D.R. *Lexrank: Graph-based Lexical Centrality as Salience in Text Summarization.* Journal of artificial intelligence research, Volume 22, pp. 457-479 (2004).
3. Horn, R.A., Johnson, C.R. *Matrix Analysis.* 575pp., Cambridge University Press, Cambridge (1990).
4. Juditsky, A., Polyak, B. *Robust Eigenvector of a Stochastic Matrix with Application to PageRank.* 51st IEEE Conference on Decision and Control, Maui, Hawaii, USA (2012).
5. Timonina-Farkas, A., Seifert, R. W. *Information Retrieval under Network Uncertainty: Robust Internet Ranking.* Operations Research, 71(6), 2328-2351.

## 6 Appendix

**Lemma 2** *Conic optimization problem*

$$f(x) = \max_{\substack{z \in \mathbb{R}^N \\ ||z||_1 \leq \varepsilon \\ |z_j| \leq \varepsilon_j}} z^T x \tag{14}$$

*is equivalent to the following minimization problem:*

$$f(x) = \min_{\lambda + \mu = x} \left\{ \varepsilon \|\lambda\|_\infty + \sum_{j=1}^{N} \varepsilon_j |\mu_j| \right\}. \tag{15}$$

*Proof* Let us dualize the optimization problem (14). First of all, notice that

$$\max_{\substack{z \in \mathbb{R}^N \\ ||z||_1 \leq \varepsilon \\ |z_j| \leq \varepsilon_j}} z^T x \iff \max_{\substack{z \in \mathbb{R}^N, t \in \mathbb{R}^N_{\{0,+\}} \\ \sum_{i=1}^{N} t_i \leq \varepsilon \\ z_j \leq t_j \\ -z_j \leq t_j \\ z_j \leq \varepsilon_j \\ -z_j \leq \varepsilon_j}} z^T x.$$

Therefore, the Lagrangian $\mathscr{L}$ can be written in the following form for dual variables $\alpha, \beta, \gamma, \eta, \nu$, where $\alpha \in \mathbb{R}_{\{0,+\}}$ and $\beta, \gamma, \eta, \nu \in \mathbb{R}^N_{\{0,+\}}$:

$$\begin{aligned} \mathscr{L} = \ & z^T x - \alpha \left( \sum_{i=1}^{N} t_i - \varepsilon \right) - \sum_{j=1}^{N} \left( \beta_j (z_j - t_j) + \gamma_j (-z_j - t_j) \right) - \\ & - \sum_{j=1}^{N} \left( \eta_j (z_j - \varepsilon_j) + \nu_j (-z_j - \varepsilon_j) \right) = \alpha \varepsilon + \sum_{j=1}^{N} (\eta_j + \nu_j) \varepsilon_j + \\ & + \sum_{j=1}^{N} z_j (x_j - \beta_j + \gamma_j - \eta_j + \nu_j) + \sum_{j=1}^{N} t_j (\beta_j + \gamma_j - \alpha). \end{aligned}$$

By strong duality, the following holds

$$\max_{\substack{z \in \mathbb{R}^N \\ ||z||_1 \leq \varepsilon \\ |z_j| \leq \varepsilon_j}} z^T x = \max_{\substack{z \in \mathbb{R}^N \\ t \in \mathbb{R}^N_{\{0,+\}}}} \min_{\substack{\beta, \gamma, \eta, \nu \in \mathbb{R}^N_{\{0,+\}} \\ \alpha \in \mathbb{R}_{\{0,+\}}}} \mathscr{L} = \min_{\substack{\beta, \gamma, \eta, \nu \in \mathbb{R}^N_{\{0,+\}} \\ \alpha \in \mathbb{R}_{\{0,+\}}}} \max_{\substack{z \in \mathbb{R}^N \\ t \in \mathbb{R}^N_{\{0,+\}}}} \mathscr{L},$$

where the following is true at the point of maximum over $z, t$:

$$x_j - \beta_j + \gamma_j - \eta_j + \nu_j = 0, \ \forall j = 1, ..., N.$$
$$\beta_j + \gamma_j - \alpha \leq 0, \ \forall j = 1, ..., N.$$

Substituting these equations into the Lagrangian and maximizing over $z$ and $t$, we get

$$\mathscr{L} = \alpha \varepsilon + \sum_{j=1}^{N} (\eta_j + \nu_j) \varepsilon_j. \tag{16}$$

Now, let us make the following change of variables:

$$\lambda_j = \beta_j - \gamma_j, \ \forall j = 1, ..., N,$$
$$\mu_j = \eta_j - \nu_j \ \forall j = 1, ..., N.$$

Notice, that $x_j = \lambda_j + \mu_j, \ \forall j = 1, ..., N$. At the point of minimum over $\alpha$, $\eta$, $\nu$ the term $\eta_j + \nu_j$ behaves as $|\mu_j|$. This happens, because at optimality $\mu_j = \eta_j$, $\nu_j = 0$ if $\mu_j \geq 0$ and $\mu_j = -\nu_j$, $\eta_j = 0$ if $\mu_j \leq 0$. Similarly, $\beta_j + \gamma_j, \ \forall j = 1, ..., N$ behaves as $|\lambda_j|, \ \forall j = 1, ..., N$ at optimality, which leads to $\alpha = \|\lambda\|_\infty$.

Hence, equation (16) under the proposed change of variables applies the statement of the Lemma 2.

**Lemma 3** *Conic optimization problem*

$$f(x) = \max_{\substack{z \in \mathbb{R}^N \\ \|z\|_2 \leq \varepsilon \\ |z_j| \leq \varepsilon_j}} z^T x \tag{17}$$

*is equivalent to the following minimization problem:*

$$f(x) = \min_{\lambda + \mu = x} \left\{ \varepsilon \|\lambda\|_2 + \sum_{j=1}^N \varepsilon_j |\mu_j| \right\}. \tag{18}$$

*Proof* Let us dualize the optimization problem (17). First of all, notice that

$$\max_{\substack{z \in \mathbb{R}^N \\ \|z\|_2 \leq \varepsilon \\ |z_j| \leq \varepsilon_j}} z^T x \iff \max_{\substack{z \in \mathbb{R}^N \\ \sqrt{\sum_{i=1}^N z_i^2} \leq \varepsilon \\ z_j \leq \varepsilon_j \\ -z_j \leq \varepsilon_j}} z^T x.$$

Therefore, the Lagrangian $\mathscr{L}$ can be written in the following form for dual variables $\alpha, \beta, \gamma$, where $\alpha \in \mathbb{R}_{\{0,+\}}$ and $\beta, \gamma \in \mathbb{R}_{\{0,+\}}^N$:

$$\mathscr{L} = z^T x - \alpha \left( \sqrt{\sum_{i=1}^N z_i^2} - \varepsilon \right) - \sum_{j=1}^N \beta_j (z_j - \varepsilon_j) - \sum_{j=1}^N \gamma_j (-z_j - \varepsilon_j).$$

By strong duality, the following holds

$$\max_{\substack{z \in \mathbb{R}^N \\ \|z\|_2 \leq \varepsilon \\ |z_j| \leq \varepsilon_j}} z^T x = \max_{z \in \mathbb{R}^N} \min_{\substack{\alpha \in \mathbb{R}_{\{0,+\}} \\ \beta, \gamma \in \mathbb{R}_{\{0,+\}}^N}} \mathscr{L} = \min_{\substack{\alpha \in \mathbb{R}_{\{0,+\}} \\ \beta, \gamma \in \mathbb{R}_{\{0,+\}}^N}} \max_{z \in \mathbb{R}^N} \mathscr{L},$$

where the following must be true at the point of maximum over $z$:

$$\frac{\partial L(z, \alpha, \beta, \gamma)}{\partial z_j} = x_j - \beta_j + \gamma_j - \alpha \frac{z_j}{\|z\|_2} = 0, \ \forall j = 1, ..., N.$$

Substituting this equation into the Lagrangian, we get

$$\mathscr{L}(z, \alpha, \beta, \gamma) = \alpha \varepsilon + \sum_{j=1}^N (\beta_j + \gamma_j) \varepsilon_j. \tag{19}$$

Now, let us make the following change of variables:

$$\lambda_j = \alpha \frac{z_j}{\|z\|_2}, \ \forall j = 1, ..., N,$$
$$\mu_j = \beta_j - \gamma_j, \ \forall j = 1, ..., N.$$

Notice, that $\alpha = \|\lambda\|_2$ and that at the point of minimum over $\alpha$, $\beta$, $\gamma$ the term $\beta_j + \gamma_j$ behaves as $|\mu_j|$. This happens, because at optimality $\beta_j = \mu_j$, $\gamma_j = 0$ if $\mu_j \geq 0$ and $\beta_j = 0$, $\gamma_j = -\mu_j$ if $\mu_j \leq 0$.

Hence, equation (19) applies the statement of the Lemma 3 under the proposed change of variables.

**Lemma 4** *(see Theorem 3.1 in [1]) For $a_i \in \mathbb{R}^{n_i}$, $\forall i = 0, ..., N$, $\xi_j \in \mathbb{R}^{n_0 \times n_j}$, $j = 1, ..., N$ the following holds:*

$$\max_{\substack{\|\xi_1\|_F \leq \varepsilon^{(\xi_1)} \\ \|\xi_2\|_F \leq \varepsilon^{(\xi_2)} \\ \dots \\ \|\xi_N\|_F \leq \varepsilon^{(\xi_N)}}} \left\| a_0 + \sum_{i=1}^N \xi_i a_i \right\|_2 = \|a_0\|_2 + \sum_{i=1}^N \varepsilon^{(\xi_i)} \|a_i\|_2.$$

*Proof*

$$\|a_0 + \sum_{i=1}^N \xi_i a_i\|_2^2 = \left( a_0 + \sum_{i=1}^N \xi_i a_i \right)^T \left( a_0 + \sum_{i=1}^N \xi_i a_i \right) =$$
$$= \|a_0\|_2^2 + \sum_{i=1}^N \|\xi_i a_i\|_2^2 + 2 \sum_{i=1}^N \|a_0^T \xi_j a_j\|_2 + 2 \sum_{i=1}^N \sum_{j=i+1}^N \|a_i^T \xi_i^T \xi_j a_j\|_2 \leq$$
$$\leq \|a_0\|_2^2 + \sum_{i=1}^N \left( \varepsilon^{(\xi_i)} \right)^2 \|a_i\|_2^2 + 2 \|a_0\|_2 \sum_{j=1}^N \varepsilon^{(\xi_j)} \|a_j\|_2 +$$
$$+ 2 \sum_{i=1}^N \sum_{j=i+1}^N \varepsilon^{(\xi_i)} \varepsilon^{(\xi_j)} \|a_i\|_2 \|a_j\|_2 = \left( \|a_0\|_2 + \sum_{i=1}^N \varepsilon^{(\xi_i)} \|a_i\|_2 \right)^2.$$

Hence,

$$\|a_0 + \sum_{i=1}^N \xi_i a_i\|_2 \leq \|a_0\|_2 + \sum_{i=1}^N \varepsilon^{(\xi_i)} \|a_i\|_2.$$

Equality holds if $\xi_i = \xi_i^* = \frac{\varepsilon^{(\xi_i)} a_0 a_i^T}{\|a_0\|_2 \|a_i\|_2}$ for $a_0 \neq 0$ (for $a_0 = 0$ one can take arbitrary $a_0: \|a_0\|_2 = 1$).

## 6.1 Proof of the Proposition 1

*Proof* Let $u = \begin{pmatrix} u_1 \\ u_2 \end{pmatrix}$, where $u_1$ is a vector of the length $N$ and $u_2$ is a vector of the length $M$. Notice, that the following equality holds due to the duality of the $l_1$-norm:

$$
||Qx - x||_1 = \left\| \begin{pmatrix} P + \xi & \zeta \\ \psi & \chi \end{pmatrix} \begin{pmatrix} x^{(1)} \\ x^{(2)} \end{pmatrix} - \begin{pmatrix} x^{(1)} \\ x^{(2)} \end{pmatrix} \right\|_1 = \left\| \begin{pmatrix} (P + \xi - I_N)x^{(1)} + \zeta x^{(2)} \\ \psi x^{(1)} + (\chi - I_M)x^{(2)} \end{pmatrix} \right\|_1 =
$$

$$
= \max_{\substack{u \in \mathbb{R}^{N+M} \\ ||u||_\infty \le 1}} \begin{pmatrix} u_1 \\ u_2 \end{pmatrix}^T \begin{pmatrix} (P + \xi - I_N)x^{(1)} + \zeta x^{(2)} \\ \psi x^{(1)} + (\chi - I_M)x^{(2)} \end{pmatrix} =
$$

$$
= \max_{\substack{u \in \mathbb{R}^{N+M} \\ ||u||_\infty \le 1}} \left( u_1^T(P - I_N)x^{(1)} + (u_1^T \xi + u_2^T \psi)x^{(1)} + (u_1^T \zeta + u_2^T(\chi - I_M))x^{(2)} \right),
$$

where $I_N$ and $I_M$ are identity matrices of sizes $N \times N$ and $M \times M$ correspondingly. Further, the function $||Qx - x||_1$ can be bounded from above in line with the triangle inequality and norm duality:

$$
||Qx - x||_1 \le ||Px^{(1)} - x^{(1)}||_1 + \max_{\substack{u \in \mathbb{R}^{N+M} \\ ||u||_\infty \le 1}} \left( u_1^T \xi + u_2^T \psi \right)x^{(1)} +
$$

$$
+ \max_{\substack{u \in \mathbb{R}^{N+M} \\ ||u||_\infty \le 1}} \left( u_1^T \zeta + u_2^T(\chi - I_M) \right)x^{(2)}.
$$

This leads to the statement:

$$
\max_{\substack{Q \in \Xi^{(l_1)}, \\ M \in \mathbb{N}}} ||Qx - x||_1 \le ||Px^{(1)} - x^{(1)}||_1 + \max_{\substack{||\xi_{|j}||_1 \le \varepsilon_j^{(\xi)} \\ \sum_{i,j} |\xi_{ij}| \le \varepsilon^{(\xi)} \\ ||\psi_{|j}||_1 \le \varepsilon_j^{(\psi)} \\ \sum_{i,j} |\psi_{ij}| \le \varepsilon^{(\psi)} \\ M \in \mathbb{N}}} \max_{\substack{u \in \mathbb{R}^{N+M} \\ ||u||_\infty \le 1}} \left( u_1^T \xi + u_2^T \psi \right)x^{(1)} +
$$

$$
+ \max_{\substack{||\zeta_{|j}||_1 \le \varepsilon_j^{(\zeta)} \\ \sum_{i,j} |\zeta_{ij}| \le \varepsilon^{(\zeta)} \\ ||\chi_{|j}||_1 \le \varepsilon_j^{(\chi)} \\ \sum_{i,j} |\chi_{ij}| \le \varepsilon^{(\chi)} \\ M \in \mathbb{N}}} \max_{\substack{u \in \mathbb{R}^{N+M} \\ ||u||_\infty \le 1}} \left( u_1^T \zeta + u_2^T(\chi - I_M) \right)x^{(2)}. \tag{20}
$$

Now, consider the following subproblem:

$$
g_1(x) = \max_{\substack{||\xi_{|j}||_1 \le \varepsilon_j^{(\xi)} \\ \sum_{i,j} |\xi_{ij}| \le \varepsilon^{(\xi)} \\ ||\psi_{|j}||_1 \le \varepsilon_j^{(\psi)} \\ \sum_{i,j} |\psi_{ij}| \le \varepsilon^{(\psi)} \\ M \in \mathbb{N}}} \max_{\substack{u \in \mathbb{R}^{N+M} \\ ||u||_\infty \le 1}} \left( u_1^T \xi + u_2^T \psi \right)x^{(1)}.
$$

Let $z = \xi^T u_1 + \psi^T u_2$. Based on the conditions of the problem (6), we can write

$$|z_j| = \left| \sum_{i=1}^{N} \xi_{ij} u_{1i} + \sum_{i=1}^{M} \psi_{ij} u_{2i} \right| \leq \sum_{i=1}^{N} |\xi_{ij}| + \sum_{i=1}^{M} |\psi_{ij}| \leq \varepsilon_j^{(\xi)} + \varepsilon_j^{(\psi)},$$

$$\|z\|_1 = \sum_{j=1}^{N} \left| \sum_{i=1}^{N} \xi_{ij} u_{1i} + \sum_{i=1}^{M} \psi_{ij} u_{2i} \right| \leq \sum_{j=1}^{N} \sum_{i=1}^{N} |\xi_{ij}| + \sum_{j=1}^{N} \sum_{i=1}^{M} |\psi_{ij}| \leq \varepsilon^{(\xi)} + \varepsilon^{(\psi)},$$

where $u_{1i}, \; \forall i = 1, ..., N$ and $u_{2i}, \; \forall i = 1, ..., M$ are $i$-th coordinates of vectors $u_1$ and $u_2$ correspondingly.

Hence, independently of $M$, we get

$$g_1(x) \leq \max_{\substack{z \in \mathbb{R}^N \\ \|z\|_1 \leq \varepsilon^{(\xi)} + \varepsilon^{(\psi)} \\ |z_j| \leq \varepsilon_j^{(\xi)} + \varepsilon_j^{(\psi)}}} z^T x^{(1)} = (\varepsilon^{(\xi)} + \varepsilon^{(\psi)}) \max_{\substack{z \in \mathbb{R}^N \\ \|z\|_1 \leq 1 \\ |z_j| \leq \frac{\varepsilon_j^{(\xi)} + \varepsilon_j^{(\psi)}}{\varepsilon^{(\xi)} + \varepsilon^{(\psi)}}}} z^T x^{(1)},$$

which is a strictly feasible conic optimization problem. Dualizing the constraints we come to

$$g_1(x) \leq (\varepsilon^{(\xi)} + \varepsilon^{(\psi)}) \min_{\lambda + \mu = x^{(1)}} \left\{ \|\lambda\|_\infty + \sum_{j=1}^{N} \frac{\varepsilon_j^{(\xi)} + \varepsilon_j^{(\psi)}}{\varepsilon^{(\xi)} + \varepsilon^{(\psi)}} |\mu_j| \right\} =$$
$$= (\varepsilon^{(\xi)} + \varepsilon^{(\psi)}) \|x^{(1)}\|_{(a)}, \tag{21}$$

where we denote the underlying norm by $\|x^{(1)}\|_{(a)}$.

Analogically, consider the other subproblem:

$$g_2(x) = \max_{\substack{\|[\zeta_j]_j\|_1 \leq \varepsilon_j^{(\zeta)} \\ \sum_{i,j} |\zeta_{ij}| \leq \varepsilon^{(\zeta)} \\ \|[\chi]_j\|_1 \leq \varepsilon^{(\chi)} \\ \sum_{i,j} |\chi_{ij}| \leq \varepsilon^{(\chi)} \\ M \in \mathbb{N}}} \max_{\substack{u \in \mathbb{R}^{N+M} \\ \|u\|_\infty \leq 1}} \left( u_1^T \zeta + u_2^T (\chi - I_M) \right) x^{(2)}.$$

Let $z = \zeta^T u_1 + (\chi^T - I_M) u_2$. Based on the conditions of the problem (6), we can write

$$|z_j| = \left| \sum_{i=1}^{N} \zeta_{ij} u_{1i} + \sum_{i=1}^{M} \chi_{ij} u_{2i} - u_{2j} \right| \leq \sum_{i=1}^{N} |\zeta_{ij}| + \sum_{i=1}^{M} |\chi_{ij}| + 1 \leq \varepsilon_j^{(\zeta)} + \varepsilon_j^{(\chi)} + 1,$$

$$\|z\|_1 = \sum_{j=1}^{M} \left| \sum_{i=1}^{N} \zeta_{ij} u_{1i} + \sum_{i=1}^{M} \chi_{ij} u_{2i} - u_{2j} \right| \leq \sum_{j=1}^{M} \sum_{i=1}^{N} |\zeta_{ij}| + \sum_{j=1}^{M} \sum_{i=1}^{M} |\chi_{ij}| + M \leq$$
$$\leq \varepsilon^{(\zeta)} + \varepsilon^{(\chi)} + M.$$

Hence,

$$g_2(x) \leq \max_{\substack{z \in \mathbb{R}^M, M \in \mathbb{N} \\ \|z\|_1 \leq \varepsilon^{(\zeta)} + \varepsilon^{(\chi)} + M \\ |z_j| \leq \varepsilon_j^{(\zeta)} + \varepsilon_j^{(\chi)} + 1}} z^T x^{(2)}.$$

Therefore, the statement of the Proposition 1 follows.

## 6.2 Proof of the Proposition 2

*Proof* Let $u = \begin{pmatrix} u_1 \\ u_2 \end{pmatrix}$, where $u_1$ is a vector of the length $N$ and $u_2$ is a vector of the length $M$.

For the case of $l_1$-norm, the following equality holds

$$||Qx - x||_1 = \left\| \begin{pmatrix} P + \xi & \zeta \\ \psi & \chi \end{pmatrix} \begin{pmatrix} x^{(1)} \\ x^{(2)} \end{pmatrix} - \begin{pmatrix} x^{(1)} \\ x^{(2)} \end{pmatrix} \right\|_1 = \left\| \begin{pmatrix} (P + \xi - I_N)x^{(1)} + \zeta x^{(2)} \\ \psi x^{(1)} + (\chi - I_M)x^{(2)} \end{pmatrix} \right\|_1 =$$
$$= \|(P + \xi - I_N)x^{(1)} + \zeta x^{(2)}\|_1 + \|\psi x^{(1)} + (\chi - I_M)x^{(2)}\|_1,$$

where $I_N$ and $I_M$ are identity matrices of the size $N \times N$ and $M \times M$ correspondingly. Using the norm duality, i.e.

$$\|(P + \xi - I_N)x^{(1)} + \zeta x^{(2)}\|_1 = \max_{\substack{u_1 \in \mathbb{R}^N \\ \|u_1\|_\infty \leq 1}} u_1^T \left( (P + \xi - I_N)x^{(1)} + \zeta x^{(2)} \right)$$

$$\|\psi x^{(1)} + (\chi - I_M)x^{(2)}\|_1 \quad = \max_{\substack{u_2 \in \mathbb{R}^M \\ \|u_2\|_\infty \leq 1}} u_2^T \left( \psi x^{(1)} + (\chi - I_M)x^{(2)} \right),$$

we choose such feasible $u_1 = u_1^*$, that $\|(P + \xi)x^{(1)} - x^{(1)}\|_1 = (u_1^*)^T(P + \xi - I_N)x^{(1)}$, $u_1^* \in \mathbb{R}^N$ and $\|u_1^*\|_\infty \leq 1$ and we fix $u_2 = \mathbb{1}_M$, where $\mathbb{1}_M$ is the $M \times 1$ vector of all-ones.
By this, we compute the lower bound for the norm $||Qx - x||_1$:

$$||Qx - x||_1 \geq (u_1^*)^T \left( (P + \xi - I_N)x^{(1)} + \zeta x^{(2)} \right) + \mathbb{1}_M^T \left( \psi x^{(1)} + (\chi - I_M)x^{(2)} \right) =$$
$$= \|(P + \xi)x^{(1)} - x^{(1)}\|_1 + (u_1^*)^T \zeta x^{(2)} + \mathbb{1}_M^T \psi x^{(1)} + \mathbb{1}_M^T \chi x^{(2)} - \mathbb{1}_M^T x^{(2)} =$$
$$= \|(P + \xi)x^{(1)} - x^{(1)}\|_1 + \mathbb{1}_M^T \psi x^{(1)} + (u_1^* - \mathbb{1}_N)^T \zeta x^{(2)},$$

where the final equation holds due to equalities (22) and (23) with $\mathbb{1}_M^T \psi x^{(1)} \geq 0$ and $(u_1^* - \mathbb{1}_N)^T \zeta x^{(2)} \leq 0$:

$$\mathbb{1}_M^T \psi = -\mathbb{1}_N^T \xi; \tag{22}$$
$$\mathbb{1}_M^T \chi = \mathbb{1}_M^T - \mathbb{1}_N^T \zeta. \tag{23}$$

Notice, that equalities (22) and (23) hold due to the column-stochasticity of the matrix $Q$ (i.e. due to conditions (2) and (3)).
Therefore, we compute the lower bound for the function $\phi_1(x)$, i.e.

$$\phi_1(x) \geq \|(P + \xi)x^{(1)} - x^{(1)}\|_1 + \mathbb{1}_M^T \psi x^{(1)} + (u_1^* - \mathbb{1}_N)^T \zeta x^{(2)}, \ \forall \xi, \psi, \zeta \in \Xi^{(l_1)}$$

As soon as this bound holds for all $\xi, \psi, \zeta$ in the perturbation set, we can set $\zeta = 0$ and $\psi = 0$ without any loss of generality. Moreover, by setting $\psi = 0$ we guarantee that conditions of A. Juditsky and B. Polyak in the work [4] are satisfied, i.e. that the matrix $P + \xi$ is column-stochastic.
Therefore, we guarantee that the $l_1$-norm formulation of A. Juditsky and B. Polyak in the work [4] is a lower bound for our optimization problem, i.e., the bound (7) follows.

### 6.3 Exact solution of the problem for new sentences

Consider the optimization problem

$$\min_{y^{(2)} \in \Sigma_M} \|y^{(2)}\|_{(b)} = \min_{\substack{y^{(2)} \in \Sigma_M \\ \lambda + \mu = y^{(2)}}} \left\{ \|\lambda\|_\infty + \sum_{j=1}^M \frac{\varepsilon_j^{(\zeta)} + \varepsilon_j^{(\chi)} + 1}{\varepsilon^{(\zeta)} + \varepsilon^{(\chi)} + M} |\mu_j| \right\}.$$

One can rewrite it as the following linear optimization problem:

$$\min_{y^{(2)} \in \Sigma_M} \|y^{(2)}\|_{(b)} = \min_{\substack{y^{(2)} \ge 0 \\ \Sigma_{j=1}^M y_j^{(2)} = 1 \\ \lambda + \mu = y^{(2)} \\ \lambda_j \le u, \forall j \\ -\lambda_j \le u, \forall j \\ \mu_j \le v_j, \forall j \\ -\mu_j \le v_j, \forall j}} \left\{ u + \sum_{j=1}^M \frac{\varepsilon_j^{(\zeta)} + \varepsilon_j^{(\chi)} + 1}{\varepsilon^{(\zeta)} + \varepsilon^{(\chi)} + M} v_j \right\}. \tag{24}$$

Let us denote $c_j = \frac{\varepsilon_j^{(\zeta)} + \varepsilon_j^{(\chi)} + 1}{\varepsilon^{(\zeta)} + \varepsilon^{(\chi)} + M}$ and solve the optimization problem (24) by its dualization.

The Lagrangian $\mathscr{L}$ for the problem (24) is

$$\begin{aligned}
\mathscr{L} &= u + \sum_{j=1}^M c_j v_j + \alpha \left( \sum_{j=1}^M y_j^{(2)} - 1 \right) + \\
&+ \sum_{j=1}^M \left( \beta_j (\lambda_j + \mu_j - y_j^{(2)}) + \gamma_j (\mu_j - v_j) + \eta_j (-\mu_j - v_j) \right) + \\
&+ \sum_{j=1}^M \left( \kappa_j (\lambda_j - u) + \nu_j (-\lambda_j - u) \right),
\end{aligned}$$

where dual variables are $\alpha \in \mathbb{R}$, $\beta \in \mathbb{R}^M$, $\gamma, \eta, \kappa, \nu \in \mathbb{R}^M_{\{0,+\}}$.

Differently, one can rewrite the Lagrangian in the following form:

$$\begin{aligned}
\mathscr{L} &= u \left( 1 - \sum_{j=1}^M \kappa_j - \sum_{j=1}^M \nu_j \right) + \sum_{j=1}^M v_j (c_j - \gamma_j - \eta_j) + \sum_{j=1}^M y_j^{(2)} (\alpha - \beta_j) + \\
&+ \sum_{j=1}^M \lambda_j (\beta_j + \kappa_j - \nu_j) + \sum_{j=1}^M \mu_j (\beta_j + \gamma_j - \eta_j) - \alpha.
\end{aligned}$$

By strong duality, the following holds

$$\min_{y^{(2)} \in \Sigma_M} \|y^{(2)}\|_{(b)} = \min_{\substack{y^{(2)} \in \mathbb{R}^M_{\{0,+\}} \\ u \in \mathbb{R}_{\{0,+\}} \\ v \in \mathbb{R}^M_{\{0,+\}} \\ \lambda, \mu \in \mathbb{R}^M}} \max_{\substack{\alpha \in \mathbb{R} \\ \beta \in \mathbb{R}^M \\ \gamma, \eta, \kappa, \nu \in \mathbb{R}^M_{\{0,+\}}}} \mathscr{L} = \max_{\substack{\alpha \in \mathbb{R} \\ \beta \in \mathbb{R}^M \\ \gamma, \eta, \kappa, \nu \in \mathbb{R}^M_{\{0,+\}}}} \min_{\substack{y^{(2)} \in \mathbb{R}^M_{\{0,+\}} \\ u \in \mathbb{R}_{\{0,+\}} \\ v \in \mathbb{R}^M_{\{0,+\}} \\ \lambda, \mu \in \mathbb{R}^M}} \mathscr{L},$$

where the following must hold true:

$$\begin{cases} 1 - \sum_{j=1}^{M} \kappa_j - \sum_{j=1}^{M} \nu_j \geq 0, \\ c_j - \gamma_j - \eta_j \geq 0, \ \forall j = 1,...,N, \\ \alpha - \beta_j \geq 0, \ \forall j = 1,...,N, \\ \beta_j + \kappa_j - \nu_j = 0, \ \forall j = 1,...,N, \\ \beta_j + \gamma_j - \eta_j = 0, \ \forall j = 1,...,N. \end{cases}$$

Therefore,

$$\min_{y^{(2)} \in \Sigma_M} \|y^{(2)}\|_{(b)} = \max_{\substack{\tilde{a} \in \mathbb{R} \\ \tilde{b} \in \mathbb{R}^M \\ \gamma, \eta, \kappa, \nu \in \mathbb{R}_{\{0,+\}}^M}} \left\{ \tilde{a}, \text{ subject to} \right.$$

$$\sum_{j=1}^{M} (\kappa_j + \nu_j) \leq 1,$$
$$c_j \geq \eta_j + \gamma_j, \ \forall j = 1,...,N,$$
$$\tilde{a} \leq \tilde{b}_j, \ \forall j = 1,...,N,$$
$$\tilde{b}_j = \kappa_j - \nu_j, \ \forall j = 1,...,N,$$
$$\left. \tilde{b}_j = \gamma_j - \eta_j, \ \forall j = 1,...,N \right\},$$

where we denote $\tilde{a} = -\alpha$ and $\tilde{b}_j = -\beta_j$. At optimality, $\nu_j = \eta_j = 0$ and, therefore, one can rewrite the optimization problem in the following equivalent form:

$$\min_{y^{(2)} \in \Sigma_M} \|y^{(2)}\|_{(b)} = \max_{\substack{\tilde{a} \in \mathbb{R} \\ \tilde{b} \in \mathbb{R}^M}} \left\{ \tilde{a}, \text{ subject to} \right.$$

$$\sum_{j=1}^{M} \tilde{b}_j \leq 1,$$
$$\tilde{a} \leq \tilde{b}_j, \ \forall j = 1,...,N,$$
$$\left. c_j \geq \tilde{b}_j, \ \forall j = 1,...,N \right\}.$$

According to this optimization problem, if $\exists \, \tilde{b}_j > \frac{1}{M}$, then $\tilde{a} < \frac{1}{M} \ \forall j$. Hence, the best possible choice for $\tilde{b}_j$ would be $\frac{1}{M}$. However, there is an additional constraint $c_j \geq \tilde{b}_j, \ \forall j = 1,...,M$. Therefore, if $\min_j \{c_j\} < \frac{1}{M}$, then $\tilde{a} = \min_j \{c_j\} < \frac{1}{M}$. One can summarize it as follows:

$$\min_{y^{(2)} \in \Sigma_M} \|y^{(2)}\|_{(b)} = \begin{cases} \frac{1}{M}, \text{ if } c_j \geq \frac{1}{M}, \ \forall j = 1,...,M \\ \min_j \{c_j\}, \text{ if } \min_j \{c_j\} < \frac{1}{M}. \end{cases}$$